\documentclass[reqno,12pt]{amsart}
\usepackage{amssymb}
\usepackage{amsmath}
\usepackage{amsfonts}

\newtheorem{theorem}{Th\'eor\`eme}[section]
\theoremstyle{plain}

\newtheorem{proposition}[theorem]{Proposition}

\numberwithin{equation}{section}
\input{tcilatex}

\begin{document}
\title{ Deux remarques sur l'espace d'interpolation}
\author{Daher Mohammad}
\address{177, Rue Gustave Courbet 77350 Le M\'{e}e Sur Seine-France}
\email{daher.mohammad@ymail.com}

\begin{abstract}
Dans ce travail, on montre que $(M(\mathbb{T}),c_{0}(\mathbb{Z}))_{\theta
}=(L^{1},c_{0}(\mathbb{Z}))_{\theta },$ $0<\theta <1.$ Dans la suite on
montre pour le couple d'interpolation $(C_{0},C_{1})$ trouv\'{e} par
Garling-Smith qu'il existe un isomorphisme $U_{\theta
}:(C_{0},C_{0}+C_{1})_{\theta ,p}\rightarrow (C_{1},C_{0}+C_{1})_{\theta ,p}$
(resp. $U_{\theta }:(C_{0},C_{0}+C_{1})_{\theta }\rightarrow
(C_{1},C_{0}+C_{1})_{\theta })$ tel que sa restriction \`{a} $C_{\theta ,p}$
(resp. \`{a} $C_{\theta })$ est un isomorphisme : $C_{\theta ,p}\rightarrow
C_{1-\theta ,p}$ (resp. $C_{\theta }\rightarrow C_{1-\theta }).$

Abstract. In this work we show that $(M(\mathbb{T}),c_{0}(\mathbb{Z}%
))_{\theta }=(L^{1},c_{0}(\mathbb{Z}))_{\theta },$ $0<\theta <1.$ In the
following we show for the interpolation couple found by Garling-Smith that
there exists an isomorphism $U_{\theta }:(C_{0},C_{0}+C_{1})_{\theta
,p}\rightarrow (C_{1},C_{0}+C_{1})_{\theta ,p}$ (resp. $U_{\theta
}:(C_{0},C_{0}+C_{1})_{\theta }\rightarrow (C_{1},C_{0}+C_{1})_{\theta })$
such that its restriction to $C_{\theta ,p}$ (resp. \`{a} $C_{\theta })$ is
an isomorphism : $C_{\theta ,p}\rightarrow C_{1-\theta ,p}$ (resp. $%
C_{\theta }\rightarrow C_{1-\theta })$
\end{abstract}

\maketitle

\bigskip AMS Clasification:45B70,46B22,46B28

\ \ \ \ \ \ \ \ \ \ \ \ \ \ \ \ \ \ \ \ \ \ \ \ \ \ \ \ \ \ \ \ \ \ \ \ \ \
\ \ \ \ \ \ \ \ \ \ \ \ \ \ \ \ \ \ \ \ \ \ \ \ \ \ \ \ \ \ \ \ \ \ \ \ \ \
\ \ \ \ \ \ \ \ \ \ \ \ \ \ \ \ \ \ \ \ \ \ \ \ \ \ \ \ \ \ \ \ \ \ \ \ \ \
\ \ \ \ \ \ \ \ \ \ \ \ \ \ \ \ \ \ \ \ \ \ \ \ \ \ \ \ \ \ \ \ \ \ \ \ \ \
\ \ \ \ \ \ \ \ \ \ \ \ \ \ \ \ \ \ \ \ \ \ \ \ \ \ \ \ \ \ \ \ \ \ \ \ \ \
\ \ \ \ \ \ \ \ \ \ \ \ \ \ \ \ \ \ \ \ \ \ \ \ \ \ \ \ \ \ \ \ \ \ \ \ \ \
\ \ \ \ \ \ \ \ \ \ \ \ \ \ \ \ \ \ \ \ \ \ \ \ \ \ \ \ \ \ \ \ \ \ \ \ \ \
\ \ \ \ \ \ \ \ \ \ \ \ \ \ \ \ \ \ \ \ \ \ \ \ \ \ \ \ \ \ \ \ \ \ \ \ \ \
\ \ \ \ \ \ \ \ \ \ \ \ \ \ \ \ \ \ \ \ \ \ \ \ \ \ \ \ \ \ \ \ \ \ \ \ \ \
\ \ \ \ \ \ \ \ \ \ \ \ \ \ \ \ \ \ \ \ \ \ \ \ \ \ \ \ \ \ \ \ \ \ \ \ \ \
\ \ \ \ \ \ \ \ \ \ \ \ \ \ \ \ \ \ \ \ \ \ \ \ \ \ \ \ \ \ \ \ \ \ \ \ \ \
\ \ \ \ \ \ \ \ \ \ \ \ \ \ \ \ \ \ \ \ \ \ \ \ \ \ \ \ \ \ \ \ \ \ \ \ \ \
\ \ \ \ \ \ \ \ \ \ \ \ \ \ \ \ \ \ \ \ \ \ \ \ \ \ \ \ \ \ \ \ \ \ \ \ \ \
\ \ \ \ \ \ \ \ \ \ \ \ \ \ \ \ \ \ \ \ \ \ \ \ \ \ \ \ \ \ \ \ \ \ \ \ \ \
\ \ \ \ \ \ \ \ \ \ \ \ \ \ \ \ \ \ \ \ \ \ \ \ \ \ \ \ \ \ \ \ \ \ \ \ \ \
\ \ \ \ \ \ \ \ \ \ \ \ \ \ \ \ \ \ \ \ \ \ \ \ \ \ \ \ \ \ \ \ \ \ \ \ \ \
\ \ \ \ \ \ \ \ \ \ \ \ \ \ \ \ \ \ \ \ \ \ \ \ \ \ \ \ \ \ \ \ \ \ \ \ \ \
\ \ \ \ \ \ \ \ \ \ \ \ \ \ \ \ \ \ \ \ \ \ \ \ \ \ \ \ \ \ \ \ \ \ \ \ \ \
\ \ \ \ \ \ \ \ \ \ \ \ \ \ \ \ \ \ \ \ \ \ \ \ \ \ \ \ \ \ \ \ \ \ \ \ \ \
\ \ \ \ \ \ \ \ \ \ \ \ \ \ \ \ \ \ \ \ \ \ \ \ \ \ \ \ \ \ \ \ \ \ \ \ \ \
\ \ \ \ \ \ \ \ \ \ \ \ \ \ \ \ \ \ \ \ \ \ \ \ \ \ \ \ \ \ \ \ \ \ \ \ \ \
\ \ \ \ \ \ \ \ \ \ \ \ \ \ \ \ \ \ \ \ \ \ \ \ \ \ \ \ \ \ \ \ \ \ \ \ \ \
\ \ \ \ \ \ \ \ \ \ \ \ \ \ \ \ \ \ \ \ \ \ \ \ \ \ \ \ \ \ \ \ \ \ \ \ \ \
\ \ \ \ \ \ \ \ \ \ \ \ \ \ \ \ \ \ \ \ \ \ \ \ \ \ \ \ \ \ \ \ \ \ \ \ \ \
\ \ \ \ \ \ \ \ \ \ \ \ \ \ \ \ \ \ \ \ \ \ \ \ \ \ \ \ \ \ \ \ \ \ \ \ \ \
\ \ \ \ \ \ \ \ \ \ \ \ \ \ \ \ \ \ \ \ \ \ \ \ \ \ \ \ \ \ \ \ \ \ \ \ \ \
\ \ \ \ \ \ \ \ \ \ \ \ \ \ \ \ \ \ \ \ \ \ \ \ \ \ \ \ \ \ \ \ \ \ \ \ \ \
\ \ \ \ \ \ \ \ \ \ \ \ \ \ \ \ \ \ \ \ \ \ \ \ \ \ \ \ \ \ \ \ \ \ \ \ \ \
\ \ \ \ \ \ \ \ \ \ \ \ \ \ \ \ \ \ \ \ \ \ \ \ \ \ \ \ \ \ \ \ \ \ \ \ \ \
\ \ \ \ \ \ \ \ \ \ \ \ \ \ \ \ \ \ \ \ \ \ \ \ \ \ \ \ \ \ \ \ \ \ \ \ \ \
\ \ \ \ \ \ \ \ \ \ \ \ \ \ \ \ \ \ \ \ \ \ \ \ \ \ \ \ \ \ \ \ \ \ \ \ \ \
\ \ \ \ \ \ \ \ \ \ \ \ \ \ \ \ \ \ \ \ \ \ \ \ \ \ \ \ \ \ \ \ \ \ \ \ \ \
\ \ \ \ \ \ \ \ \ \ \ \ \ \ \ \ \ \ \ \ \ \ \ \ \ \ \ \ \ \ \ \ \ \ \ \ \ \
\ \ \ \ \ \ \ \ \ \ \ \ \ \ \ \ \ \ \ \ \ \ \ \ \ \ \ \ \ \ \ \ \ \ \ \ \ \
\ \ \ \ \ \ \ \ \ \ \ \ \ \ \ \ \ \ \ \ \ \ \ \ \ \ \ \ \ \ \ \ \ \ \ \ \ \
\ \ \ \ \ \ \ \ \ \ \ \ \ \ \ \ \ \ \ \ \ \ \ \ \ \ \ \ \ \ \ \ \ \ \ \ \ \
\ \ \ \ \ \ \ \ \ \ \ \ \ \ \ \ \ \ \ \ \ \ \ \ \ \ \ \ \ \ \ \ \ \ \ \ \ \
\ \ \ \ \ \ \ \ \ \ \ \ \ \ \ \ \ \ \ \ \ \ \ \ \ \ \ \ \ \ \ \ \ \ \ \ \ \
\ \ \ \ \ \ \ \ \ \ \ \ \ \ \ \ \ \ \ \ \ \ \ \ \ \ \ \ \ \ \ \ \ \ \ \ \ \
\ \ \ \ \ \ \ \ \ \ \ \ \ \ \ \ \ \ \ \ \ \ \ \ \ \ \ \ \ \ \ \ \ \ \ \ \ \
\ \ \ \ \ \ \ \ \ \ \ \ \ \ \ \ \ \ \ \ \ \ \ \ \ \ \ \ \ \ \ \ \ \ \ \ \ \
\ \ \ \ \ \ \ \ \ \ \ \ \ \ \ \ \ \ \ \ \ \ \ \ \ \ \ \ \ \ \ \ Mots cl\'{e}%
s: Interpolation des espaces $\mathbf{h}^{p}$ et $L^{p}$

\begin{center}
\textsc{Introduction}
\end{center}

D\'{e}sigons par $M(\mathbb{T})$ l'espace des mesure sur $\mathbb{T}.$ Pour $%
\mu \in M(\mathbb{T})$ et $t\in \mathbb{T}$ on d\'{e}finit $\mu _{t}$ par $%
\left\langle f,\mu _{t}\right\rangle =\left\langle f_{t},\mu \right\rangle ,$
$f\in C(\mathbb{T}),$ ici $f_{t}(y)=f(y-t),$ $y\in \mathbb{T}.$

Rappelons qu'il existe une suite $(K_{n})_{n\geq 0}$ born\'{e}e dans $L^{1}(%
\mathbb{T})$ telle que, pour pour tout Banach $X$ et toute $f\in C(\mathbb{T}%
,X),$ $(f\ast K_{n})_{n\geq 0}$ converge vers $f$ dans $C(\mathbb{T},X).$

D\'{e}signons par $E$ l'espace des mesures $\mu $ d\'{e}finies sur $\mathbb{T%
}$ \`{a} valeurs complexes telle que $\widehat{\mu }(n)\rightarrow 0$ quand $%
\left\vert n\right\vert \rightarrow +\infty .$

Pour les d\'{e}finitions des espaces d'interpolation $A_{\theta },A_{\theta
,p}$ nous r\'{e}f\'{e}rons \`{a} \cite{Ber-Lof}.

Pour tout $\theta \in \left] 0,1\right[ $ notons $X_{\theta }=(L^{1}(\mathbb{%
T}),c_{0}(\mathbb{Z}))_{\theta },$ $Y_{\theta }=(M(\mathbb{T}),\ell ^{\infty
}(\mathbb{Z}))_{\theta }$ et $Z_{\theta }=(E,c_{0}(\mathbb{Z}))_{\theta }$ .

\textbf{\ }Soient$\ B_{0}=L^{1}(\mathbb{T)},$ $B_{0}=\widehat{c_{0}(\mathbb{Z%
})}$ l'espace des transform\'{e}es de Fourier de $L^{1}(\mathbb{T)},$ muni
de la norme de $c_{0}(\mathbb{Z}).$ D'apr\`{e}s \cite{Blas-Xu}, $X_{\theta }$
contient $c_{0}$ isomorphiquement.

\begin{proposition}
\label{ccv}Pour tout $\theta \in \left] 0,1\right[ $ on a $X_{\theta
}=Z_{\theta }$ isom\'{e}triquement.
\end{proposition}

Consid\'{e}rons $\mu \in E$ et $t\in \mathbb{T}.$ Comme $\mu \in Y_{\theta }$%
, $\mu _{t}\in Y_{\theta }.$ Montrons que l'application $U$ : $t\in \mathbb{T%
}\rightarrow \mu _{t}\in Y_{\theta }$ est continue. Soit $\varepsilon >0.$
Il existe $n_{0}\in \mathbb{N}$ tel que si $\left\vert n\right\vert \geq
n_{0}$, $\left\vert \widehat{\mu }(n)\right\vert <\varepsilon /2.$

Remarquons que $\widehat{\mu }_{t}(n)=e^{int}$ $\widehat{\mu }(n),$ $n\in 
\mathbb{Z},t\in \mathbb{T},$ donc pour tout $t,t^{\prime }\in \mathbb{T}$%
\begin{eqnarray}
\sup_{n\in \mathbb{Z}}\left\vert \widehat{\mu }_{t}(n)-\widehat{\mu }%
_{t^{\prime }}(n)\right\vert &=&\sup_{n\in \mathbb{Z}}\left\vert \widehat{%
\mu }(n)\right\vert \left\vert e^{int}-e^{int^{\prime }}\right\vert  \notag
\\
&=&\max (\sup_{\left\vert n\right\vert <n_{0}}\left\vert \widehat{\mu }%
(n)\right\vert \left\vert e^{int}-e^{int^{\prime }}\right\vert
,\sup_{\left\vert n\right\vert \geq n_{0}}\left\vert \widehat{\mu }%
(n)\right\vert \left\vert e^{int}-e^{int^{\prime }}\right\vert )  \notag \\
&\leq &\max (\sup_{\left\vert n\right\vert <n_{0}}\left\vert \widehat{\mu }%
(n)\right\vert \left\vert e^{int}-e^{int^{\prime }}\right\vert ,\varepsilon
).  \label{fff}
\end{eqnarray}

Il en r\'{e}sulte que l'application $t\in \mathbb{T}\rightarrow (\widehat{%
\mu }_{t}(n))_{n\in \mathbb{Z}}\in c_{0}(\mathbb{Z})$ est continue.

D'autre part, d'apr\`{e}s (\ref{fff}) $\left\Vert \mu _{t}-\mu _{t^{\prime
}}\right\Vert _{Z_{\theta }}\leq C_{\theta }\left\Vert \mu _{t}-\mu
_{t^{\prime }}\right\Vert _{M(\mathbb{T})}^{1-\theta }(\sup_{n\in \mathbb{Z}%
}\left\vert \widehat{\mu }_{t}(n)-\widehat{\mu }_{t^{\prime }}(n)\right\vert
)^{\theta }$ ($C_{\theta }$ est une constante)$.$ D'apr\`{e}s ce qui pr\'{e}c%
\`{e}de l'application $U$ est continue. Par cons\'{e}quent $K_{m}\ast
U(0)\rightarrow _{m\rightarrow +\infty }\mu $ dans $Y.$ Comme pour tout $m,$ 
$K_{m}\ast U(0)\in L^{1}(\mathbb{T})\subset X_{\theta }$ et $X_{\theta }$
est un spus-espace isom\'{e}trique de $Y_{\theta }$ d'apr\`{e}s \cite[Lemme
3.8]{Da1}, alors $\mu \in X_{\theta }.\blacksquare $

\begin{theorem}
\label{tr}\cite[th.2]{Gar-Smi}Il existe un couple d'interpolation $%
(C_{0},C_{1})$ tel que $C_{j}$ est isomorphe \`{a} $\ell ^{1},$ $j=0,1$ et $%
(C_{0},C_{1})_{\theta },$ $(C_{0},C_{1})_{\theta ,p}$ contiennent $c_{0}$
isomorphiquement, pour tout $\theta \in \left] 0,1\right[ $ et tout $p\in %
\left[ 1,+\infty \right[ .$
\end{theorem}

\begin{theorem}
\label{yi}Il existe un isomorphisme $U_{\theta }:(C_{0},C_{0}+C_{1})_{\theta
,p}\rightarrow (C_{1},C_{0}+C_{1})_{\theta ,p}$ (resp. $U_{\theta
}:(C_{0},C_{0}+C_{1})_{\theta }\rightarrow (C_{1},C_{0}+C_{1})_{\theta })$
tel que sa restriction \`{a} $C_{\theta ,p}$ (resp. \`{a} $C_{\theta })$ est
un isomorphisme : $C_{\theta ,p}\rightarrow C_{1-\theta ,p}$ (resp. $%
C_{\theta }\rightarrow C_{1-\theta }).$
\end{theorem}

D\'{e}monstration.

Montrons que $(C_{0},C_{0}+C_{1})_{\theta ,p}$ est isomorphe \`{a} $%
(C_{1},C_{0}+C_{1})_{\theta ,p}.$

L'ensemble des \'{e}l\'{e}ments dans $c_{0}$ sous la forme ($\pm 1,\pm
1,...,\pm 1,0,0,0,...)$ est d\'{e}nombrable. Consid\'{e}rons $(r_{n})_{n\geq
1}$ une num\'{e}ration de cet ensemble. Notons $\varepsilon
_{n}=(1+\left\Vert r_{n}\right\Vert _{\ell ^{1}})^{-n},$ $n\geq 1.$

D'apr\`{e}s \cite[th.2]{Gar-Smi}

\begin{equation}
C_{0}=\left\{ (\underset{n}{\dsum }(a_{n}e_{n}+b_{n}r_{n}),\underset{n}{%
\dsum }\varepsilon _{n}b_{n}e_{n},\underset{n}{\dsum }c_{n}e_{n});\text{ }%
(a_{n})_{n\geq 1},(b_{n})_{n\geq 1},(c_{n})_{n\geq 1}\in \ell ^{1}\right\} ,
\label{ll}
\end{equation}%
\begin{equation}
C_{1}=\left\{ (\underset{n}{\dsum }(a_{n}e_{n}+c_{n}r_{n}),\underset{n}{%
\dsum }b_{n}e_{n},\underset{n}{\dsum }\varepsilon _{n}c_{n}e_{n});\text{ }%
(a_{n})_{n\geq 1},(b_{n})_{n\geq 1},(c_{n})_{n\geq 1}\in \ell ^{1}\right\} .
\label{yy}
\end{equation}

Soit $U_{0}:$ $C_{0}\rightarrow C_{1},$ l'isomorphisme d\'{e}finie par $%
U_{0}(x)=(\underset{n}{\dsum }(a_{n}e_{n}+b_{n}r_{n}),\underset{}{\underset{n%
}{\dsum }c_{n}e_{n},\underset{n}{\dsum }}\varepsilon _{n}b_{n}e_{n}),$ $x=(%
\underset{n}{\dsum }(a_{n}e_{n}+b_{n}r_{n}),\underset{n}{\dsum }\varepsilon
_{n}b_{n}e_{n},\underset{n}{\dsum }c_{n}e_{n})\in C_{0}$.

Soit d'autre part, $U_{1}:C_{0}+C_{1}\rightarrow C_{0}+C_{1}$ d\'{e}finie
par $U_{1}(x+y)=U_{0}(x)+(U_{0})^{-1}(y),$ $x\in C_{0},y\in C_{1}.$ Montrons
que $U_{1}$ est une applicationet et injective. Pour cela; soient $x=(%
\underset{}{\underset{n}{\dsum }}(a_{n}e_{n}+b_{n}r_{n}),\underset{n}{\dsum }%
\varepsilon _{n}b_{n}e_{n},\underset{n}{\dsum }c_{n}e_{n},),u=(\underset{n}{%
\dsum }(a_{n}^{\prime }e_{n}+b_{n}^{\prime }r_{n}),\underset{n}{\dsum }%
\varepsilon _{n}b_{n}^{\prime }e_{n},\underset{n}{\dsum }c_{n}^{\prime
}e_{n})\in C_{0}$, $y=(\underset{n}{\dsum }(\alpha _{n}e_{n}+\gamma
_{n}r_{n}),\underset{n}{\dsum }\beta _{n}e_{n},\underset{n}{\dsum }%
\varepsilon _{n}\gamma _{n}e_{n}),v=(\underset{n}{\dsum }(\alpha
_{n}^{\prime }e_{n}+\gamma _{n}^{\prime }r_{n}),\underset{n}{\dsum }\beta
_{n}^{\prime }e_{n},\underset{n}{\dsum }\varepsilon _{n}\gamma _{n}^{\prime
}e_{n})\in C_{1}.$ Remarquons que $x+y=u+v$ signifie que%
\begin{eqnarray*}
&&\underset{n}{(\dsum }\left[ (a_{n}+\alpha _{n})e_{n}+(b_{n}+\gamma
_{n})r_{n}\right] ,\underset{n}{\dsum }(\varepsilon _{n}b_{n}+\beta
_{n})e_{n},\underset{n}{\dsum }(c_{n}+\varepsilon _{n}\gamma _{n})e_{n}) \\
&=&\underset{n}{(\dsum }\left[ (a_{n}^{\prime }+\alpha _{n}^{\prime
})e_{n}+(b_{n}^{\prime }+\gamma _{n}^{\prime }r_{n})\right] ,\underset{n}{%
\dsum }(\varepsilon _{n}b_{n}^{\prime }+\beta _{n}^{\prime })e_{n},\underset{%
n}{\dsum }(c_{n}^{\prime }+\varepsilon _{n}\gamma _{n}^{\prime })e_{n}).
\end{eqnarray*}

Il en r\'{e}sulte que 
\begin{eqnarray}
\underset{n}{\dsum }\left[ (a_{n}+\alpha _{n})e_{n}+(b_{n}+\gamma _{n}r_{n})%
\right] &=&\underset{n}{\dsum }\left[ (a_{n}^{\prime }+\alpha _{n}^{\prime
})e_{n}+(b_{n}^{\prime }+\gamma _{n}^{\prime }r_{n})\right] ,  \notag \\
\underset{n}{\dsum }(\varepsilon _{n}b_{n}+\beta _{n})e_{n} &=&\underset{n}{%
\dsum }(\varepsilon _{n}b_{n}^{\prime }+\beta _{n}^{\prime })e_{n},
\label{ii} \\
\underset{n}{\dsum }(c_{n}+\varepsilon _{n}\gamma _{n})e_{n} &=&\underset{n}{%
\dsum }(c_{n}^{\prime }+\varepsilon _{n}\gamma _{n}^{\prime })e_{n}.  \notag
\end{eqnarray}

D'autre part, ($U_{0})^{-1}(y)=(\underset{n}{\dsum }(\alpha _{n}e_{n}+\gamma
_{n}r_{n}),\underset{n}{\dsum \varepsilon _{n}}\gamma _{n}e_{n},\underset{n}{%
\dsum }\beta _{n}e_{n}),$ ($U_{0})^{-1}(v)=(\underset{n}{\dsum }(\alpha
_{n}^{\prime }e_{n}+\gamma _{n}^{\prime }r_{n}),\underset{n}{\dsum
\varepsilon _{n}}\gamma _{n}^{\prime }e_{n},\underset{n}{\dsum }\beta
_{n}^{\prime }e_{n}),$ donc 
\begin{eqnarray}
&&U_{0}(x)+(U_{0})^{-1}(y)  \notag \\
&=&(\underset{n}{\dsum }\left[ a_{n}+\alpha _{n})e_{n}+(b_{n}+\gamma
_{n})r_{n}\right] ,\underset{}{\underset{n}{\dsum }(c_{n}+\varepsilon
_{n}\gamma _{n})e_{n},\underset{n}{\dsum }}(\varepsilon _{n}b_{n}+\beta
_{n})e_{n}),  \notag \\
&&U_{0}(u)+(U_{0})^{-1}(v)  \label{jj} \\
&=&\underset{n}{(\dsum }\left[ a_{n}^{\prime }+\alpha _{n}^{\prime
})e_{n}+(b_{n}^{\prime }+\gamma _{n}^{\prime })r_{n}\right] ,\underset{}{%
\underset{n}{\dsum }(c_{n}^{\prime }+\varepsilon _{n}\gamma _{n}^{\prime
})e_{n},\underset{n}{\dsum }}(\varepsilon _{n}b_{n}^{\prime }+\beta
_{n}^{\prime })e_{n}).  \notag
\end{eqnarray}

Il est clair d'apr\`{e}s (\ref{ii}) et (\ref{jj}) que $%
U_{0}(x)+(U_{0})^{-1}(y)=U(u)+(U_{0})^{-1}(v)$ si et seulement si $x+y=u+v,$
c'est-\`{a}-dire que $U_{0}$ est une application et injective.

Comme $U_{0}:C_{0}\rightarrow C_{1},$ $U_{1}$ $:C_{0}+C_{1}\rightarrow
C_{0}+C_{1}$ sont isomorphismes, d'apr\`{e}s \cite[th.3.1.2]{Ber-Lof} $%
U_{\theta }:(C_{0},C_{0}+C_{1})_{\theta ,p}\rightarrow
(C_{1},C_{1}+C_{0})_{\theta ,p}$ est un isomrphisme$.$

Remarquons que la restriction de $U_{1}$ \`{a} $C_{1}$ est un isomorphisme : 
$C_{1}\rightarrow C_{0}$, par cons\'{e}quent la restriction de $U_{\theta }$ 
\`{a} $C_{\theta ,p}$ est un isomorphisme : $C_{\theta ,p}\rightarrow
(C_{1},C_{0})_{\theta ,p}=(C_{1},C_{0})_{1-\theta ,p}.$

Par un argument analogue on montre que $U_{\theta
}:(C_{0},C_{0}+C_{1})_{\theta }\rightarrow (C_{1},C_{0}+C_{1})_{\theta }$
est un isomorphisme et la restriction de $U_{\theta }$ \`{a} $C_{\theta }$
est un isomorphisme : $C_{\theta }\rightarrow C_{1-\theta }.$

\end{document}